# A few equalities involving integrals of the logarithm of the Riemann $\varsigma$-function and equivalent to the Riemann hypothesis III. Exponential weight functions.

Sergey K. Sekatskii, Stefano Beltraminelli, and Danilo Merlini

**Abstract.** This paper is a continuation of our recent papers with the same title, arXiv:0806.1596v1 [math.NT], arXiv:0904.1277v1 where a number of integral equalities involving integrals of the logarithm of the Riemann $\varsigma$-function were introduced and it was shown that some of them are equivalent to the Riemann hypothesis. A few new equalities of this type, which this time involve exponential functions, are established, and for the first time we have found equalities involving the integrals of the logarithm of the Riemann $\varsigma$-function taken exclusively along the real axis.

  Some of the equalities we have found are tested numerically. In particular, an integral equality involving the logarithm of $|\varsigma(1/2+it)|$ and a weight function $\cosh^{-1}(\pi t)$ is shown numerically to be correct up to the 80 digits. For exponential weight function $e^{-at}$, the possible contribution of the Riemann function zeroes non-lying on the critical line is rigorously estimated and shown to be extremely small, in particular, smaller than trillion of digits, $10^{-10^{13}}$, for $a = 4\pi$.

# 1. Introduction

In recent papers [1, 2] we analyzed integrals $\int_{b-i\infty}^{b+i\infty} g(z)\ln(\varsigma(z))dz$ and have established a number of equalities of the type $\int_{b-i\infty}^{b+i\infty} g(z)\ln(\varsigma(z))dz = f(b)$ which were proven to be equivalent to the Riemann Hypothesis (RH; here $\varsigma(z)$ is the Riemann zeta-function, see e.g. [3] for definitions and discussion of the general properties of this function). Our paper [2] has been finished just with the recipe how different other relations of this type can be established. In this Note we analyze a few additional equalities equivalent to the RH. These new equalities employ exponential weight functions and, in our opinion, are especially interesting. The main part of the results contained in [1, 2] as well as the first part of the content of the present article, is currently published in *Ukrainian Math. J.* [4].

## 2. Integral equalities equivalent to the Riemann hypothesis. Hyperbolical sine and cosine.

Let us consider the rectangular contour $C$ with the vertices $b-iX$, $b+iX$, $b+X+iX$, $b+X-iX$ with real $b$>-2 (with this choice we avoid unnecessary complications with the trivial Riemann zeroes) and real $X \to \infty$, i.e. the same as in [1, 2], introduce the function $g(z) = \dfrac{i}{\cos(a(z-b))}$, $a$ is real positive, and consider the contour integral $\int_C \ln(\varsigma(z))g(z)dz$. For $z-b = x+iy$ where $x, y$ are real, we have $\cos(a(z-b)) = \cos(ax)\cosh(ay) - i\sin(ax)\sinh(ay)$ which means that for large $y$ $|\cos(a(z-b))|^{-1} = O(e^{-a|y|})$ provided $\arg(z-b) \neq 0$ and $\arg(z-b) \neq \pi$. This behaviour, together with the known asymptotic of $\ln(\varsigma(z))$, guaranties the disappearance of the integral taken along the external lines of the contour (i.e. along lines other than its left border): the problems might appear only for real positive *z-b*, but for such a case we have $\ln(\varsigma(z)) \cong 2^{-z}$, so it is enough to avoid such values of *X* for which *g(z)* has poles when $X \to \infty$.

In the interior of the contour we have simple poles at the points $z = b + \pi/(2a) + \pi n/a$ where *n* is an integer or zero, and, if *b<1*, also a simple pole of the Riemann function at *z=1*. If *b<1,* in the interior of the contour we also can have a number of zeroes of the Riemann function, and we definitely have an infinite number of them if *b<1/2*. The contribution of the poles of *g(z)* to the contour integral value is, by a residue theorem, equal to $2\pi i(i\sum_{n=0}^{\infty}\frac{(-1)^{n+1}}{a}\ln(\varsigma(b+\pi/(2a)+\pi n/a))$. To find the corresponding contribution of the pole of the Riemann function we apply the generalized Littlewood theorem [1, 2, 4, 5]: for $b<1$ it is equal to $2\pi i\int_{b}^{1}\frac{i}{\cos(a(z-b))}dz=-2\pi\int_{0}^{1-b}\frac{1}{\cos(ax)}dx=\frac{\pi}{a}\ln\frac{1-\sin(a(1-b))}{1+\sin(a(1-b))}$; we take $a(1-b)<\pi/2$ to avoid the problems with this integral. Analogously, the order $l_k$ zero of the Riemann function $\rho=\sigma_k+it_k$, for $b<\sigma_k$, contributes

$$-2\pi i\int_{b+it_k}^{\sigma_k+it_k}\frac{il_k}{\cos(a(z-b))}dz=2\pi\int_{0}^{\sigma_k-b}\frac{l_k}{\cos(ax)\cosh(at_k)-i\sin(ax)\sinh(at_k)}dx$$

to the contour integral value. (Of course, the relation $\int\frac{dx}{\cos x}=\frac{1}{2}\ln\frac{1+\sin x}{1-\sin x}$ can be also used to estimate these integrals but apparently this is not of a much help). Pairing the complex conjugate zeroes we see that due to the symmetry of the distribution of the Riemann function zeroes the imaginary part of these contributions vanishes and thus, collecting everything together, we have the following equality: for $-2\le b<1$

$$\int_{0}^{\infty}\frac{\ln|\varsigma(b+it)|}{\cosh(at)}dt=\pi\sum_{\rho,t_k>0}2l_k\int_{0}^{\sigma_k-b}\frac{\cos(ax)\cosh(at_k)}{\cos^2(ax)\cosh^2(at_k)+\sin^2(ax)\sinh^2(at_k)}dx$$
$$+\frac{\pi}{2a}\ln\frac{1-\sin(a(1-b))}{1+\sin(a(1-b))}+\pi\sum_{n=0}^{\infty}\frac{(-1)^n}{a}\ln(\varsigma(b+\pi/(2a)+\pi n/a)) \qquad (1).$$

The sum here is taken over all zeroes $\rho=\sigma_k+it_k$ with $b<\sigma_k$ taken into account their multiplicities; when obtaining (1) we paired the contributions of the complex conjugate zeroes $\rho=\sigma_k\pm it_k$ hence $t_k>0$. Throughout the paper, if there is a zero or pole of the Riemann function lying on the integration path for an integral of the type $\int_{b-i\infty}^{b+i\infty}\ln(\varsigma(b+it))g(z)dz$ or similar, corresponding integral is understood as an appropriate limit. This does not create any problem (as it is sometimes said, a logarithmical peculiarity is "integrable" because of the finiteness of the limit $\lim_{x\to 0}\int_{x}^{a}\ln z dz$) and is quite common in the field, cf. [1, 2, 4, 6 - 8].

Of course, for $b\ge 1$ we have simply

$$\frac{1}{\pi}\int_{0}^{\infty}\frac{\ln|\varsigma(b+it)|}{\cosh(at)}dt=\frac{1}{a}\sum_{n=0}^{\infty}(-1)^n\ln(\varsigma(b+\pi/(2a)+\pi n/a)) \qquad (2).$$

This is easy to see that if we take $a(1-b) < \pi/2$ than the integrand in
$\int_0^{\sigma_k - b} \frac{\cos(ax)\cosh(at_k)}{\cos^2(ax)\cosh^2(at_k) + \sin^2(ax)\sinh^2(at_k)} dx$ is always positive (the value $\sigma_k - b$ can not exceed *1-b*) and thus we have proven the following

THEOREM 1. An equality
$$\frac{a}{\pi}\int_0^\infty \frac{\ln|\varsigma(b+it)|}{\cosh(at)} dt = \frac{1}{2}\ln\frac{1-\sin(a(1-b))}{1+\sin(a(1-b))} + \sum_{n=0}^\infty (-1)^n \ln(\varsigma(b + \pi/(2a) + \pi n/a)) \quad (3)$$
where *b*, *a* are real positive numbers such that $1 > b \geq 1/2$, $a(1-b) < \pi/2$ holds true for some *b* if and only if there are no Riemann function zeroes with $\sigma > b$. For *b=1/2* this equality is equivalent to the Riemann hypothesis.

It is interesting to take *b=1/2* and consider the limit $a \to \pi$. There are no Riemann function zeroes with Re*s*=1 hence the positivity of the contributions of Riemann function zeroes non-lying on the critical line is still certain and it remains only to consider the limit $\lim_{a \to \pi}(\frac{1}{2}\ln\frac{1-\sin(a(1-b))}{1+\sin(a(1-b))} + \ln(\varsigma(1/2 + \frac{\pi}{2a})))$. This can be done without problems (the first term tends to $\ln(\delta/4)$ and second to $\ln(2\pi/\delta)$ where $\delta = \pi - a$ small) and we obtain the following equality equivalent to the Riemann hypothesis
$$\int_0^\infty \frac{\ln|\varsigma(1/2+it)|}{\cosh(\pi t)} dt = \ln\frac{\pi}{2} + \sum_{n=2}^\infty (-1)^{n+1} \ln(\varsigma(n)) \quad (4).$$
Here those logarithms of the Riemann function which correspond to even *n* can be expressed via Bernoulli numbers $B_{2m}$ ($m = 1,2,3...$): $\varsigma(2m) = (-1)^{m+1}(2\pi)^{2m}\frac{B_{2m}}{2(2m)!}$, and the following remarkable property of the Riemann $\varsigma$-function $\prod_{n=2}^\infty \varsigma(n) = C = 2.29485...$, that is $\sum_{n=2}^\infty \ln(\varsigma(n)) = \ln C = 0.8306...$ can be also used for calculations, see P. 16 of [9]. Here *C* is the residue of the pole at *s=1* of the Dirichlet series whose coefficients *g(n)* are the numbers of non-isomorphic Abel's group of order *n*. The results of the numerical testing of eq. (4) is presented in the Section 5 of the paper.

For completeness it is desirable to establish similar criteria pertinent for integrals over the argument of the Riemann function rather than logarithm of its module. For this not even but odd function *g(z)* is to be considered, cf. [2]. Unfortunately, one is not able to take $g(z) = \dfrac{i}{\sin(a(z-b))}$ because the resulting integral over the line $b+it$ diverges at $t=0$. Situation can be corrected introducing, similarly to Wang [6] (see also discussion in [1]), a deformed contour which avoids the point $z=0$, but this does not seem very interesting to us. Thus here we limit ourselves with the simplest suitable odd function $g(z) = \dfrac{-i(z-b)}{\cos^2(a(z-b))}$, $a$ is real positive, and consider the same contour as above. In the interior of the contour we have *second order* poles at the points $z = b + \pi/(2a) + \pi n/a$ where $n$ is an integer or zero. The contribution of these poles to the contour integral value is equal to

$$\dfrac{2\pi i}{a^2}(-i\sum_{n=0}^{\infty}(\ln(\varsigma(b+\pi/(2a)+\pi n/a))) + (\dfrac{\pi}{2a}+\dfrac{\pi n}{a})\dfrac{\varsigma'(b+\pi/(2a)+\pi n/a)}{\varsigma(b+\pi/(2a)+\pi n/a)}).$$

The corresponding contribution of the pole of the Riemann function is equal to

$$2\pi i\int_b^1 \dfrac{-i(z-b)}{\cos^2(a(z-b))}dz = 2\pi\int_0^{1-b}\dfrac{x}{\cos^2(ax)}dx\,;$$

we take $a(1-b) < \pi/2$ to avoid the problems with this integral. Analogously, the order $l_k$ zero of the Riemann function $\rho = \sigma_k + it_k$, for $b < \sigma_k$, contributes

$$2\pi i\int_{b+it_k}^{\sigma_k+it_k}\dfrac{i(z-b)l_k}{\cos^2(a(z-b))}dz = -2\pi l_k\int_0^{\sigma_k - b}\dfrac{x+it_k}{(\cos(ax)\cosh(at_k)-i\sin(ax)\sinh(at_k))^2}dx$$

to the contour integral value. Real part of this integral, which is interesting now, is given by

$$I(t_k) = -2\pi l_k\int_0^{\sigma_k - b}\dfrac{x(\cos^2(ax)\cosh^2(at_k) - \sin^2(ax)\sinh^2(at_k)) - 2t_k\sin(ax)\cos(ax)\sinh(at_k)\cosh(at_k)}{(\cos^2(ax)\cosh^2(at_k) + \sin^2(ax)\sinh^2(at_k))^2}dx$$

Collecting everything together, we have the following equality: for $-2 \leq b < 1$

$$\int_0^{\infty}\dfrac{t\arg(\varsigma(b+it))}{\cosh^2(at)}dt = \sum_{\sigma_k > b, t_k > 0}I(t_k) + \pi\int_0^{1-b}\dfrac{x}{\cos^2(ax)}dx +$$

$$\dfrac{\pi}{a^2}\sum_{n=0}^{\infty}(\ln(\varsigma(b+\pi/(2a)+\pi n/a))) + (\dfrac{\pi}{2a}+\dfrac{\pi n}{a})\dfrac{\varsigma'(b+\pi/(2a)+\pi n/a)}{\varsigma(b+\pi/(2a)+\pi n/a)})$$

This is well known that non-trivial Riemann function zeroes not lying on the critical line, if they exist, have very large values of $|t_k|$: from numerical calculations performed by Wedeniwski (as cited by Ramaré and Saouter [10]), $|t_k| > 3.3 \cdot 10^9$. The integrand in $I(t_k)$ is always negative provided $2a(\sigma_k - b) < \pi$ and $t_k > \dfrac{x}{2\sin(ax)}\operatorname{cth}(at_k)\cos(ax)$. For large $|t_k|$ the latter is always true if one avoids the cases when *a* is extremely small; definitely, it suffices to take $a \geq 10^{-9}$ and thus we have established the following

THEOREM 1a. An equality

$$\int_0^\infty \frac{t \arg(\varsigma(b+it))}{\cosh^2(at)} dt = \pi \int_0^{1-b} \frac{x}{\cos^2(ax)} dx +$$

$$\frac{\pi}{a^2} \sum_{n=0}^\infty (\ln(\varsigma(b+\pi/(2a)+\pi n/a))) + (\frac{\pi}{2a} + \frac{\pi n}{a}) \frac{\varsigma'(b+\pi/(2a)+\pi n/a)}{\varsigma(b+\pi/(2a)+\pi n/a)})$$

where *b*, *a* are real positive numbers such that $1 > b \geq 1/2$, $a(1-b) < \pi/2$, $a \geq 10^{-9}$, holds true for some *b* if and only if there are no Riemann function zeroes with $\sigma > b$. For *b=1/2* this equality is equivalent to the Riemann hypothesis.

Similarly to the case of an integral involving the logarithm of the module of the Riemann $\varsigma$–function, we can consider the limit $a \to \pi$ which for *b=1/2* gives an interesting formula

$$\int_0^\infty \frac{t \arg(\varsigma(1/2+it))}{\cosh^2(\pi t)} dt = \frac{1}{\pi}(\ln(\pi) + \frac{\gamma}{2} - 1) + \frac{1}{\pi} \sum_{n=1}^\infty (\ln(\varsigma(1+n)) + (n+\frac{1}{2}) \frac{\varsigma'(1+n)}{\varsigma(1+n)}),$$ the

correctness of which has been tested numerically. Of course, $f(x) \equiv \ln(\varsigma(1+x)) + (x+1/2) \frac{\varsigma'(1+x)}{\varsigma(1+x)} = \frac{d}{dx}((x+1/2)\ln(\varsigma(1+x)))$ hence Euler – MacLaurin summation formula looks quite appropriate to calculate $g = \frac{1}{\pi} \sum_{n=1}^\infty (\ln(\varsigma(1+n)) + (n+\frac{1}{2}) \frac{\varsigma'(1+n)}{\varsigma(1+n)})$. Retaining only the first terms we have $g = \frac{1}{\pi}(\frac{3}{2}\ln(\varsigma(2)) + \frac{1}{2}(\ln(\varsigma(2)) + \frac{3}{2}\frac{\varsigma'(2)}{\varsigma(2)}) - \frac{B_2}{2!}\frac{df}{dx}(1) + ...)$ In this approximation we obtain *g* = -0.2994461350 while the computation of the above sum with 50 first terms gives -0.2991293527. This shows a fast convergence of the MacLaurin expansion and deserves further numerical investigation.

## 3. Integral equalities equivalent to the Riemann hypothesis. Exponents.

*3.1. Integrals along lines b+it with $1 > b \geq 1/2$*

Let us consider somewhat different contour than above, viz. the rectangular contour *C* with the vertices *b*, *b + iX*, *b + X + iX*, *b + X* with real *b>-2* (with this choice we avoid the trivial Riemann zeroes) and real $X \to \infty$, introduce the function $g(z) = i\exp(i(a(z-b)))$, *a* is real positive, and consider the contour integral $\int_C g(z) \ln(\varsigma(z)) dz$. Along the line *(b, b+iX)* we have an integral $\int_0^\infty e^{-at}\ln(\varsigma(b+it))dt$ and along the line *(b, b+X)* we have an integral

$$i\int_0^\infty (\cos(ax)+i\sin(ax))\ln(\varsigma(b+x))dx = -\int_0^\infty \sin(ax)\ln(\varsigma(b+x))dx + i\int_0^\infty \cos(ax)\ln(\varsigma(b+x))dx$$

On the border of the contour we have the pole of the Riemann function at *b=1*. If *b<1*, in the interior of the contour we also can have a number of zeroes of the Riemann function, and we definitely have an infinite number of them if *b<1/2*. The order $l_k$ zero of the Riemann function $\rho = \sigma_k + it_k$, for $b < \sigma_k$, contributes

$$-2\pi i l_k \int_0^{\sigma_k - b} e^{-at_k}(i\cos(ax)-\sin(ax))dx = \frac{2\pi l_k e^{-at_k}}{a}[\sin(a(\sigma_k - b))+i(1-\cos(a(\sigma_k - b)))] \quad \text{to}$$

the contour integral value. Hence we established the equality:

$$\int_0^\infty e^{-at}\ln(\varsigma(b+it))dt - \int_0^\infty \sin(ax)\ln(\varsigma(b+x))dx + i\int_0^\infty \cos(ax)\ln(\varsigma(b+x))dx =$$

$$\frac{2\pi}{a}\sum_{\rho, t_k > 0, \sigma_k > b} l_k e^{-at_k}[\sin(a(\sigma_k - b))+i(1-\cos(a(\sigma_k - b)))] - \quad (5).$$
$$\frac{\pi}{a}(\sin(a(1-b))+i(1-\cos(a(1-b))))$$

The value $\sigma_k - b$ can not exceed *1-b* and hence if we take $a \leq \frac{\pi}{1-b}$ both real and imaginary summation terms in (5) are always positive. For imaginary terms, $a \leq \frac{2\pi}{1-b}$ can be taken. To complete the calculation, the choice of the branch of an argument is to be specified. Selecting for real *x>1* $\arg(\varsigma(x)) = 0$ and for $1 > b \geq 1/2$ $\arg(\varsigma(b+i\varepsilon)) = -\pi$ where real positive $\varepsilon \to 0$, we establish the following

THEOREM 2. Equalities

$$\int_0^\infty e^{-at}\ln|\varsigma(b+it)|dt - \int_0^\infty \sin(ax)\ln|\varsigma(b+x)|dx + \frac{\pi}{a}\sin(a(1-b)) = 0 \qquad (6)$$

$$\int_0^\infty e^{-at}\arg(\varsigma(b+it))dt + \int_0^\infty \cos(ax)\ln|\varsigma(b+x)|dx + \frac{\pi}{a}(1-\cos(a(1-b))) = 0 \qquad (7)$$

where *b*, *a* are real positive numbers such that $1 > b \geq 1/2$, $a(1-b) \leq \pi$ for (6) and $a(1-b) \leq 2\pi$ for (7), hold true for some *b* if and only if there are no Riemann function zeroes with $\sigma > b$. For *b=1/2* these equalities are equivalent to the Riemann hypothesis.

If necessary, in (7) for *b=1/2* one can take $a = \pi$ thus eliminating any problem with the integration over real axis (pole of the Riemann function coincides with the zero of cosine):

$$\int_0^\infty e^{-\pi t} \arg(\varsigma(1/2+it))dt + \int_0^\infty \cos(\pi x) \ln|\varsigma(1/2+x)| dx + 1 = 0 \qquad (7a)$$

Similarly, for (6) the choice $a = 2\pi$ is useful:

$$\int_0^\infty e^{-2\pi t} \ln|\varsigma(1/2+it)| dt - \int_0^\infty \sin(2\pi x) \ln|\varsigma(1/2+x)| dx = 0 \qquad (6a).$$

The results of the numerical testing of eq. (6a, 7a) are presented in the Section 5 of the paper.

### 3.2. *Integrals along lines b+it with $b < 1/2$ and b>1*

Summands $1 - \cos(a(\sigma_k - b)) = 2\sin^2\left(\frac{a(\sigma_k - b)}{2}\right)$, appearing in (5) for the Theorem 2 concerning the imaginary part of the integral, i.e. for $\int_0^\infty e^{-at} \arg(\varsigma(b+it))dt$, are of course non-negative for *all* values of *a*. Unfortunately, one still can not put $a > \frac{2\pi}{1-b}$, e.g. $a > 4\pi$ for *b=1/2*, into the corresponding condition of the Theorem 2, because, if so, it is still possible that Riemann zeroes have *exactly* such exceptional values of $\sigma_k$ for which $1 - \cos(a(\sigma_k - b)) = 0$ hence eq. (7) might hold true while there *are* Riemann zeroes with $\sigma_k > 1/2$. Of course, this is a degenerate case: the smallest change of *a* would again restore an inequality.

The same non-negativity of $1 - \cos(a(\sigma_k - b))$ summands enables to establish other interesting equalities. Suppose we take *b<1/2* and selected *a* in such a manner that $1 - \cos(a(1/2-b)) = 0$, viz. $a = \frac{2\pi n}{1/2 - b}$ with *n*=1, 2,... Then the contributions of all Riemann zeroes lying on the critical line are equal to zero while those of the Riemann zeroes with $\sigma \neq 1/2$ are positive provided the function $1 - \cos(a(\sigma_k - b))$ is not equal to zero for $0 < \sigma_k < 1/2$, $1/2 < \sigma_k < 1$. The latter condition can be fulfilled if one takes $b \leq 0$ and $a = \frac{2\pi n}{1/2 - b}$ with a positive integer $n \leq 1 - 2b$ (so, of course, always $a \leq 4\pi$ and we again can not take $a > 4\pi$ to avoid the same "degenerate case" which has been discussed above). Thus we have established the following

THEOREM 3. Equality

$$\int_0^\infty e^{-at} \arg(\varsigma(b+it)) dt + \int_0^\infty \cos(ax) \ln|\varsigma(b+x)| dx + \frac{\pi}{a}(1-\cos(a(1-b))) - \frac{\pi}{a} \sum_{k=1}^{k<[-b/2]} (1-\cos(a(-2k-b))) = 0 \quad (8),$$

where $b$ is any real non-positive number and $a = \frac{2\pi n}{1/2 - b}$ where $n$ is a positive integer $n \leq 1 - 2b$, is equivalent to the Riemann hypothesis and holds true if and only if there are no Riemann function zeroes with $\sigma \neq 1/2$. In this equality we start from the following initial value of an argument function: $\arg(\varsigma(b+i\varepsilon)) = -\pi + [-b/2] \cdot \pi$, where real positive $\varepsilon \to 0$.

Here, of course, $\sum_{k=1}^{k<[-b/2]} (1-\cos(a(-2k-b)))$ is the contribution of trivial zeroes while $1-\cos(a(1-b))$ is that of a simple pole. As far as we know, this is the first integral equality involving the integrals along the line $b+it$ with $b<1/2$ which is equivalent to the Riemann hypothesis.

Up to now in this paper as well as in [1, 2, 4] we always considered integrals involving the logarithm of the Riemann $\varsigma$-function following the conditions of generalized Littlewood theorem [1, 2, 4, 5]. In particular, this means that the value of $\arg(\varsigma(b+it))$ jumps on the $2\pi l$ each time when we, moving along the line $(b, b+i\infty)$ starting from $b$, pass an ordinate of Riemann zero of the $l$-th order, i.e. if there exists zero of the $l$-th order at $\sigma + it$ with $\sigma > b$. Correspondingly, the function $\arg(\varsigma(b+it))$ is not continuous but only piece-wise continuous if such zeroes exist (quite the contrary, function $\ln|\varsigma(b+it)|$ is always continuous – of course, apart from the cases when some zeroes lie exactly on the line $(b, b+i\infty)$). If we want to re-express integrals along the line $(b, b+i\infty)$ as integrals along the line *1-b+it*, using a functional equation in the form $\varsigma(s) = \pi^{s-1/2} \frac{\Gamma(1/2 - s/2)}{\Gamma(s/2)} \varsigma(1-s)$ or $\varsigma(s) = \frac{1}{\pi}(2\pi)^s \sin\frac{\pi s}{2} \Gamma(1-s)\varsigma(1-s)$, this circumstance should be taken into account and integrals $\int_0^\infty g(t)\arg(\varsigma(b+it))dt$ should be properly modified to become integrals $\int_0^\infty g(t)\arg^*(\varsigma(b+it))dt$. Here we use an asterisk * sign to underline that the function $\arg^*(\varsigma(b+it))$ is, in the sense of Remark 1 presented below, continuous along the line $(0, +\infty)$.

Such modification is straightforward: clearly, if real functions $f^*(x)$, $g(x)$ are continuous on the segment [A, B] and function $f(x)$ coincides with $f^*(x)$ on [A, C] and is equal to $f^*(x)+h$ on (C, B] ($h$ is a constant, $A \leq C \leq B$), we have $\int_A^B g(x)f(x)dx = \int_A^B g(x)f^*(x)dx + h\int_C^B g(x)dx$, and this is trivially generalized for complex numbers and improper integrals having counting number of discontinuities. For our particular case we have, considering one discontinuity at $t=t_k$: $\int_0^\infty e^{-at}\arg(\varsigma(b+it))dt = \int_0^\infty e^{-at}\arg^*(\varsigma(b+it))dt + \frac{2\pi d_k}{a}e^{-at_k}$, and the corresponding generalization gives:

$$\int_0^\infty e^{-at}\arg(\varsigma(b+it))dt = \int_0^\infty e^{-at}\arg^*(\varsigma(b+it))dt + \sum_{\rho, t_k > 0, \sigma_k > b} \frac{2\pi d_k}{a}e^{-at_k} \qquad (9).$$

In other words, integrals involving continuous (see Remark below) branch of an argument and those requiring by the conditions of the generalized Littlewood theorem differ from each other on the value $\Sigma_a(b)$: $\Sigma_a(b) = \sum_{\rho, t_k > 0, \sigma_k > b} \frac{2\pi d_k}{a}e^{-at_k}$. If we take $b \leq 0$, it does not any more depend on $b$ being a constant depending only on the value of $a$; below we will use a notation $\Sigma_a = \sum_{\rho, t_k > 0} \frac{2\pi d_k}{a}e^{-at_k}$ for such a case.

REMARK 1. Here the following remark is at place. We do *not* have $2\pi d$ jump of the argument if some zero occurs *exactly* on the line $(b, b+i\infty)$ and, correspondingly, no summand appears in $\Sigma_a(b)$ for the case, although the argument function is not continuous but jumps on $\pi$ (or larger value for zeroes of larger order) passing through such zero. This jump is, however, quite "natural" and "locally determined" having nothing common with the "artificial nonlocal" $2\pi d$ jumps depending on the behaviour of the function far from the point of jump (viz., existence of a zero at $\sigma + it$ with $\sigma > b$). Having this in mind we will continue to speak about "continuous" argument function for these cases.

Using (9), we get instead of (5):
$$\int_0^\infty e^{-at} \ln^*(\varsigma(b+it))dt - \int_0^\infty \sin(ax)\ln(\varsigma(b+x))dx + i\int_0^\infty \cos(ax)\ln(\varsigma(b+x))dx =$$
$$\frac{2\pi}{a} \sum_{\rho, t_k>0, \sigma_k>b} l_k e^{-at_k}[\sin(a(\sigma_k-b))-i\cos(a(\sigma_k-b))] - \frac{\pi}{a}(\sin(a(1-b))+i(1-\cos(a(1-b))))$$
(5a).

Real part of the r.h.s. here remains the same while for an imaginary part we have
$$\int_0^\infty e^{-at} \arg^*(\varsigma(b+it))dt + \frac{2\pi}{a}\sum_{\rho,t_k>0,\sigma_k>b} l_k e^{-at_k}\cos(a(\sigma_k-b))$$
$$+\int_0^\infty \cos(ax)\ln|\varsigma(b+x)|dx + \frac{\pi}{a}(1-\cos(a(1-b))) - \frac{\pi}{a}\sum_{k=1}^{k\leq[-b/2]}(1-\cos(a(-2k-b))) = 0 \quad (10),$$

where again the following branch of the argument should be taken: we start from $\arg(\varsigma(b+i\varepsilon)) = -\pi + [-b/2]\cdot\pi$, where real positive $\varepsilon \to 0$, and then the argument is continuous along the line $(b, b+i\infty)$. Using above definition of $\Sigma_a(b)$, it is useful to rewrite (10) also in the following form:
$$\int_0^\infty e^{-at}\arg^*(\varsigma(b+it))dt - \frac{\pi}{a}\sum_{\rho,t_k>0,\sigma_k>b} l_k e^{-at_k}\sin^2\frac{a(\sigma_k-b)}{2} + \Sigma_a(b)$$
$$+\int_0^\infty \cos(ax)\ln|\varsigma(b+x)|dx + \frac{\pi}{a}(1-\cos(a(1-b))) - \frac{\pi}{a}\sum_{k=1}^{k\leq[-b/2]}(1-\cos(a(-2k-b))) = 0 \quad (10a).$$

As we mentioned above, integrals $\int_0^\infty e^{-at}\arg^*(\varsigma(b+it))dt$ with $b<0$ can be rewritten using an appropriate functional equation. Further, it is useful to express the appearing here integral along the line *1-b-it*, where an argument function is now continuous in the sense of Remark 1, using again the equality (7), which for the case at hand is unconditionally true: $\int_0^\infty e^{-at}\arg(\varsigma(1-b-it))dt = \int_0^\infty \cos(ax)\ln|\varsigma(1-b+x)|dx$. Then, noting that $\int_0^\infty e^{-at}\arg(\pi^{b-1/2+it})dt = \frac{\ln\pi}{a^2}$ and that, of course, $\arg(\Gamma(1/2-b/2-it/2)) = -\arg(\Gamma(1/2-b/2+it/2))$, we rewrite (8) as

$$\frac{\ln \pi}{a^2} - \int_0^\infty e^{-at}(\arg(\Gamma(1/2-b/2+it/2)) + \arg(\Gamma(b/2+it/2)))dt + \Sigma_a$$

$$+ \int_0^\infty \cos(ax)(\ln|\varsigma(b+x)| + \ln|\varsigma(1-b+x)|)dx + \frac{\pi}{a}(1-\cos(a(1-b))) - \frac{\pi}{a}\sum_{k=1}^{k<[-b/2]}(1-\cos(a(-2k-b))) = 0$$

(11).

Here the initial value of an argument of $\Gamma(1/2-b/2+i\varepsilon)$ is equal to zero while that of $\Gamma(b/2+i\varepsilon)$ is the same as in the Theorem 3.

If one selects $-b$ to be a half integer, $-b=n-1/2$, the difference of the arguments of the gamma-functions appearing in (9) is an integer and then the relation $\Gamma(x+1) = x\Gamma(x)$ can be used a few times to simplify the expressions further. As a particular example we can take, say, $b=-1/2$ to get the following equalities equivalent to the Riemann hypothesis: for $a = 2\pi$

$$\frac{\ln \pi}{4\pi^2} - \int_0^\infty e^{-2\pi t}(2\arg(\Gamma(3/4+it/2)) - \arg(-\frac{1}{4}+\frac{it}{2}))dt + \Sigma_{2\pi}$$

$$+ \int_0^\infty \cos(2\pi x)(\ln|\varsigma(-1/2+x)| + \ln|\varsigma(3/2+x)|)dx + 1/2 = 0 \qquad (12a)$$

and for $a = 4\pi$

$$\frac{\ln \pi}{16\pi^2} - \int_0^\infty e^{-4\pi t}(2\arg(\Gamma(3/4+it/2)) - \arg(-\frac{1}{4}+\frac{it}{2}))dt + \Sigma_{4\pi}$$

$$+ \int_0^\infty \cos(4\pi x)(\ln|\varsigma(-1/2+x)| + \ln|\varsigma(3/2+x)|)dx + 1/4 = 0 \qquad (12b).$$

Another interesting particular example is to take $b=0$, $a = 4\pi$ in (11): repeating the same as above we get

$$\frac{\ln \pi}{16\pi^2} - \int_0^\infty e^{-4\pi t}(\arg(\Gamma(1/2+it/2)) + \arg(\Gamma(it/2)))dt + \Sigma_{4\pi}$$

$$+ \int_0^\infty \cos(4\pi x)(\ln|\varsigma(x)| + \ln|\varsigma(1+x)|)dx = 0,$$

and here Legendre's duplication formula $\Gamma(1/2+it/2)\Gamma(it/2) = \sqrt{\pi}2^{1-it}\Gamma(it)$ can be used to establish that the equality

$$\frac{\ln(2\pi)}{16\pi^2} - \int_0^\infty e^{-4\pi t}\arg(\Gamma(it))dt + \Sigma_{4\pi} + \int_0^\infty \cos(4\pi x)(\ln|\varsigma(x)| + \ln|\varsigma(1+x)|)dx = 0 \qquad (13)$$

is equivalent to the Riemann hypothesis. As far we know, these are the first integral equalities involving the integrals containing the values of the Riemann function taken exclusively along the real axis which are equivalent to the Riemann hypothesis. Some integrals appearing in Eqs. (12, 13) are elementary or can be handled using an appropriate table, e.g. [10]. We will not consider these aspects in the present paper.

Somewhat analogous equality can be established on the basis of the relation

$$\int_0^\infty e^{-at} \ln|\varsigma(b+it)| \, dt - \int_0^\infty \sin(ax) \ln|\varsigma(b+x)| \, dx + \frac{\pi}{a} \sin(a(1-b)) =$$
$$\frac{2\pi d}{a} \sum_{\rho, t_k > 0} e^{-at_k} \sin(a(\sigma_k - b))$$

concerning the real part of the integral in Theorem 2. If we take $-2 \leq b \leq 0$ and select $a$ in such a manner that $\sin(a(1/2-b)) = 0$, i.e. $a = \frac{\pi n}{1/2 - b}$ with $n=1, 2\ldots$, all contributions of Riemann zeroes lying on the critical line are equal to zero. At the same time, the contributions of Riemann zeroes with $\sigma_k = 1/2 + \alpha$ and $\sigma_k = 1/2 - \alpha$ exactly compensate each other, and thus we have established that the equality $\int_0^\infty e^{-at} \ln|\varsigma(b+it)| \, dt - \int_0^\infty \sin(ax) \ln|\varsigma(b+x)| \, dx + \frac{\pi}{a} \sin(a(1-b)) = 0$ with $a = \frac{\pi n}{1/2 - b}$ is *unconditionally* true for $-2 \leq b \leq 0$; and this is easy to modify it for $b<-2$. Taking $0 < b \leq 1/4$ we easily establish the following

THEOREM 3a. Equality

$$\int_0^\infty e^{-at} \ln|\varsigma(b+it)| \, dt - \int_0^\infty \sin(ax) \ln|\varsigma(b+x)| \, dx + \frac{\pi}{a} \sin(a(1-b)) = 0 \qquad (14),$$

where $a, b$ are real positive numbers such that $1/4 \geq b > 0$ and $a = \frac{\pi n}{1/2 - b}$, $n$ is a positive integer such that $n \leq \frac{1}{2b} - 1$, holds true for some $b$ if and only if there are no Riemann function zeroes with $\sigma < b$. In particular, the equality $\int_0^\infty e^{-4\pi t} \ln|\varsigma(1/4+it)| \, dt - \int_0^\infty \sin(4\pi x) \ln|\varsigma(1/4+x)| \, dx = 0$ holds true if and only if there are no zeroes with $\sigma < 1/4$.

The proof is clear from above: Riemann zeroes with $1/2 \geq \sigma \geq b$ contribute nothing to the contour integral value while if $\sigma_k < b$, the contribution of its "partner" Riemann zero with $\sigma = 1 - \sigma_k$ lying in the interior of the contour is not cancelled by that of the zero with $\sigma_k$ simply because the latter lies outside the contour. The conditions of the theorem ensure that all these contributions have the same sign. Unfortunately, similar theorem can not be formulated for $1/2 > b > 1/4$ because no one value of $a = \frac{\pi n}{1/2 - b}$ can ensure the same sign of all Riemann zeroes contributions for such a case.

Finally, starting from (10), similar theorem can be established for an integral involving continuous, in the sense of Remark 1, $\arg^*(\varsigma(b+it))$ function with $0<b<1/2$. If, in conditions of (10) and considering the same contour as above, we put $a = \frac{\pi}{1-2b}$, all contributions of Riemann $\varsigma$-function zeroes in the sum $\frac{2\pi}{a} \sum_{\rho, t_k > 0} l_k e^{-at_k} \cos(a(\sigma_k - b))$ cancel out due to the circumstances that $\cos(a(1/2-b)) = 0$ and that zeroes are symmetrical with respect to the line $\sigma = 1/2$. If, however, there is some zero with $\sigma_k < b$, the contribution of its partner zero lying at $1-\sigma_k$ is not cancelled and is definitely negative, so we obtain the following

THEOREM 3b. Equality
$$\int_0^\infty e^{-at} \arg^*(\varsigma(b+it)) dt + \int_0^\infty \cos(ax) \ln|\varsigma(b+x)| dx + \frac{\pi}{a}(1 - \cos(a(1-b))) = 0 \qquad (15),$$
where $a, b$ are real positive numbers such that $1/2 > b > 0$ and $a = \frac{\pi(2n+1)}{1-2b}$, $n$ is a positive integer such that $n \leq \frac{1}{2b} - 1$, holds true for some $b$ if and only if there are no Riemann function zeroes with $\sigma < b$.

For $b<0$ similar relation is unconditionally true. As you see, in conditions of Theorem 3b the complete interesting range of $b$ values can be used: here we have not the limitation $1/4 \geq b > 0$ appearing in conditions of Theorem 3a.

## 4. Rigorous bound for the contribution of remaining Riemann zeroes not lying on the critical line

Especially simple weight function appearing in the integrals considered in the previous Section makes the problem of estimation of the maximal possible contribution of remaining Riemann function zeroes non-lying on the critical line a rather simple one. For example, for the real part we know that

$$\int_0^\infty e^{-at} \ln|\varsigma(1/2+it)|\, dt - \int_0^\infty \sin(ax) \ln|\varsigma(1/2+x)|\, dx + \frac{\pi}{a}\sin(a/2) =$$
$$\frac{2\pi}{a} \sum_{\rho, t_k > 0} l_k e^{-at_k} \sin a(\sigma_k - 1/2))$$

where the sum is over all zeroes non-lying on the critical line and having $\sigma > 1/2$. How large the r.h.s. here can be? First, we have that $\frac{2\pi}{a}|\sum_{\rho, t_k>0} l_k e^{-at_k} \sin a(\sigma_k - 1/2)| \leq \frac{2\pi}{a} \sum_{\rho, t_k>0} l_k e^{-at_k}$. This last sum is given by Stieltjes integral $I = \frac{2\pi}{a} \int_0^\infty e^{-at} dN(t)$, where $N(t)$ is a discontinuous counting function which counts the number of Riemann zeroes having the real part $\sigma_k > 1/2$ and $T \geq t_k > 0$ taking into account the order of zeroes. The integration by parts is possible here, and such an integration readily gives $I = \frac{2\pi}{a^2} \int_0^\infty e^{-at} N(t) dt$. About the function $N(t)$ we know that it is equal to zero when $t < T = 3.3 \cdot 10^9$ [11]. Further, when $t > T$ it can not exceed the function $\tilde{N}(t)/2 - \tilde{N}(T)/2$, where $\tilde{N}(t)$ is a usual zero counting function which counts zeroes in the strip $1 > \sigma_k > 0$, $T \geq t_k > 0$. ($N(t)$ is equal to $\tilde{N}(t)/2 - \tilde{N}(T)/2$ if all zeroes with $t>T$ do not lye on the critical line). Thus the contribution at question can not exceed

$$I = \frac{\pi}{a^2} \int_T^\infty e^{-at}(\tilde{N}(t) - \tilde{N}(T))\, dt \qquad (16).$$

Now for the function $\tilde{N}(t)$ we use the known formula

$$\tilde{N}(t) = \frac{t}{2\pi}(\ln \frac{t}{2\pi} - 1) + \frac{7}{8} + Q(t) \qquad (17)$$

with the estimation of the reminder given by Backlund in 1918 [12] and apparently remaining the best known one:
$$|Q(t)| < 0.137 \ln t + 0.443 \ln \ln t + 4.350 \qquad (18).$$
Then integration in (16) with (17, 18) is trivial and we obtain that

$$|I| < \frac{1}{2a^3}(2\pi \cdot 0.137 \ln T + \frac{1}{a}\ln T + 2\pi \cdot 4.35)e^{-aT} + \frac{2\pi \cdot 0.443}{2a^3}\int_T^\infty \frac{e^{-at}}{\ln^2 t} dt + (\frac{2\pi \cdot 0.137}{2a^3} + \frac{1}{2a^5})\int_T^\infty \frac{e^{-at}}{t} dt$$

Of course, $\int_T^\infty \frac{e^{-at}}{\ln^2 t} dt < \frac{1}{a}\frac{e^{-aT}}{\ln^2 T}$, $\int_T^\infty \frac{e^{-at}}{t} dt < \frac{1}{a}\frac{e^{-aT}}{T}$ and thus

$$|I| < \frac{1}{2a^3}(0.861 \ln T + \frac{1}{a}\ln T + 27.332 + \frac{2.784}{a\ln^2 T} + \frac{0.861}{aT} + \frac{1}{a^3 T})e^{-aT} \qquad (19).$$

Thus for, say, $a = 1$ we obtain $|I| < 34.1 e^{-3.3 \cdot 10^9} = \dfrac{34.1}{10^{1.433 \cdot 10^9}}$, and taking the maximal possible value of $a$, here $a = 2\pi$, we get $|I| < 0.101 e^{-2.074 \cdot 10^{10}} = \dfrac{0.101}{10^{0.9 \cdot 10^{10}}}$ i.e a precision with at least 1.43 milliard decimals in the first case and of 9 milliards of decimals in the second case.

Similar estimation can be easily obtained for imaginary part, that is for

$$\int_0^\infty e^{-at} \arg(\varsigma(1/2 + it)) dt + \int_0^\infty \cos(ax) \ln|\varsigma(1/2 + x)| dx + \frac{\pi}{a}(1 - \cos(a/2)) =$$

$\dfrac{2\pi}{a} \sum_{\rho, t_k > 0} l_k e^{-at_k}(1 - \cos(a(\sigma_k - 1/2)))$, by noting that $|1 - \cos(a(\sigma_k - 1/2))| \leq 2$ hence the numerical coefficient is simply twice larger than above. But here we can put $a = 4\pi$ thus obtaining $|I| < 0.025 e^{-4.14 \cdot 10^{10}}$, i.e a precision of 13 trillions of decimals (here 1 trillion = $10^{12}$). Reminding the remark made in the beginning of the Section 3.2, we even "almost can" take an arbitrary large value of $a$ here thus making a remaining term, which is already small enough, "almost arbitrary small".

REMARK2. There is a more recent calculation of Gourdon and Demichel where it is reported that the first $k = 10^{13}$ Riemann zeroes are located on the critical line, but we were unable to get an exact value of $T$ from the paper [13]. If one will use a known estimation $T \cong 2\pi k / \ln k$ [3], page 214, as small error term as $I < c \cdot 10^{-1.318 \cdot 10^{13}}$ follows for the sum $\dfrac{2\pi}{a} \sum_{\rho, t_k > 1} l_k e^{-at_k} \sin a(\sigma_k - 1/2))$ over Riemann zeroes non-lying on the critical line, here $c = \dfrac{1}{2a^3}(0.861 \ln T + \dfrac{1}{a} \ln T + 27.332 + \dfrac{2.784}{a \ln^2 T} + \dfrac{0.861}{aT} + \dfrac{1}{a^3 T})$, and similarly for the $\dfrac{2\pi}{a} \sum_{\rho, t_k > 1} l_k e^{-at_k}(1 - \cos(a(\sigma_k - 1/2)))$ sum.

## 5. Numerical results

Equality (4) seems well-suited for numerical treatments, and we actually obtained surprising and interesting results with a relatively modest computational effort. With less than 14 hours of CPU time spent for the evaluation of (4) on a standard Personal Computer we obtained for the difference between its l.h.s. and r.h.s. the value $1.162608547 \times 10^{-81}$. In other words, the match between l.h.s. and r.h.s. is exact until the 80th digit. The integration was accomplished using the standard procedure *Integrate* in *Mathematica* with the working precision set to 200 and the maximal level of recursion set to 350. For the variable *t* we considered the integration interval [0, 200] and we truncated the sum at *n*=1000, where *n* is the number of the Riemann function zero in the sequence of these zeroes arranging in an ascending order, i.e. such that $\rho_1 = 14.13...$, $\rho_2 = 21.02...$, and so on. For an illustration, we present the computer transcript concerning this calculation in Fig. 1.

The exponential case (equations (6) and (7)) is numerically harder to analyze. With the same requirements as before the algorithm didn't converge, so we had to reduce the required calculation's precision. Because the first integrand of (6a) and (7a) decays very fast, we can limit the integration interval of the variable *t* to [0, 50]. We need also to reduce the working precision (set to 50) and maximal recursion (set to 70). For the second integral in (6a) and (7a) we set the integration interval for the variable *x* to [0, 200], and set the working precision and maximal recursion to 50. After 1.5 hour CPU time we obtained for (6a) the result of $8.8044282 \times 10^{-32}$ and for (7a) after more than 8 hours - $1.14767406 \times 10^{-31}$. So in both cases we have a correspondence for at least the first 30 digits.

The last case that we have treated numerically is equality (8) with $a = 2\pi$ and $b = -\frac{7}{2}$, i.e. we also consider the first trivial Riemann zero. In this case, we calculated the first integral with the same parameters as before. For the second integral we extended the integration's interval to [0, 250] and the maximal recursion to 100. The result $8.74533283 \times 10^{-30}$ is similar to two previous cases. Note that here $\arg(\varsigma(-7/2 + i\varepsilon)) = 0$ and, from Theorem 3, the contribution of the pole exactly compensates that of the first trivial zero hence in fact the two integrals, when summing, compensate each other up to 30 digits as found in the computations.

## 6. Conclusions

In this paper, which is a continuation of our earlier works [1, 2, 4], we have established a number of new criteria involving the integrals of the logarithm of the Riemann $\varsigma$–function and equivalent to the Riemann hypothesis, this time with the exponential weight functions. For the case of $\cosh^{-1}(\pi t)$ weight function, the rapid decrease of the integrand and some other peculiarities of the equality enabled to achieve a very large precision, up to the 80 digits, when the equality has been tested numerically. Exponential weight functions lead to rather simple expressions for the contribution of the Riemann function zeroes not lying on the critical line to the contour integral value. This equality enabled to obtain a rigorous estimation of the possible error which was shown to be extremely small.

These results are interesting but, we believe, the most exciting is the founding of completely new possibilities to establish such integral equalities equivalent to the RH, which involve the integrals over the logarithm of the Riemann $\varsigma$-function taken along the vertical line $b+it$ lying outside the critical strip $0<b<1$; see Section 3.2. Using the functional equation for the Riemann $\varsigma$-function, the latter can be rewritten in a form involving the integrals over the logarithm of the Riemann $\varsigma$-function taken exclusively along the real axis. Our opinion is that all the possibilities open by these results are not yet fully exploited, and this case deserves additional and careful analysis.

REFERENCES

[1] S. K. Sekatskii, S. Beltraminelli, and D. Merlini, A few equalities involving integrals of the logarithm of the Riemann $\varsigma$-function and equivalent to the Riemann hypothesis, arXiv:0806.1596v1 [math.NT].

[2] S. K. Sekatskii, S. Beltraminelli, and D. Merlini, A few equalities involving integrals of the logarithm of the Riemann $\varsigma$-function and equivalent to the Riemann hypothesis II, arXiv:0904.1277v1

[3] E. C. Titchmarsh and E. R. Heath-Brown, The theory of the Riemann Zeta-function, Oxford, Clarendon Press, 1988.

[4] S. K. Sekatskii, S. Beltraminelli, and D. Merlini, On equalities involving integrals of the logarithm of the Riemann $\varsigma$-function and equivalent to the Riemann hypothesis, Ukrainian Math. J., 64, (2012), 218.

[5] E. C. Titchmarsh, The theory of functions, Oxford, Oxford Univ. Press, 1939.

[6] F. T. Wang, A note on the Riemann Zeta-function, Bull. Amer. Math. Soc. 52 (1946), 319.

[7] V. V. Volchkov, On an equality equivalent to the Riemann hypothesis, Ukrainian Math. J., 47, (1995), 422.

[8] Balazard M., Saias E. and Yor M., Notes sur la fonction de Riemann, Advances in Mathematics, 143, (1999), 284.

[9] D. B. Zagier, Zetafunktionen und quadratische Körper (Eine Einführung in die höhere Zahlentheorie), Springer-Verlag, Berlin-Heidelberg - New York ,1981.

[10] I. S. Gradshtein et I. M. Ryzhik, Tables of integrals, series and products, Academic, New York, 1990.

[11] Cited from : C. Ramaré and Y. Saouter, Short effective intervals containing primes. J. of Number Theor., 98 (2003) 10-33.

[12] R. J. Backlund, Uber die Nullstellen der Riemannschen Zetafunktion, Acta Math., 41 (1918), 345 – 375.

[13] X. Gourdon, The $10^{13}$ first zeroes of the Riemann Zeta Function, and zeroes computation at very large height (2004). http://numbers.computation.free.fr/Constants/miscellaneous/zetazeroes1e13_1e24.pdf

S. K. Sekatskii, Laboratoire de Physique de la Matière Vivante, IPSB, BSP 408, Ecole Polytechnique Fédérale de Lausanne, CH1015 Lausanne-Dorigny, Switzerland.

E-mail : serguei.sekatski@epfl.ch

S. Beltraminelli, CERFIM, Research Center for Mathematics and Physics, PO Box 1132, 6600 Locarno, Switzerland.

E-mail: Stefano.beltraminelli@ti.ch

D. Merlini, CERFIM, Research Center for Mathematics and Physics, PO Box 1132, 6600 Locarno, Switzerland.

E-mail: merlini@cerfim.ch

```
g = Log[Abs[Zeta[1/2 + I*t]]] / Cosh[Pi*t]
```

$$\text{Log}\left[\text{Abs}\left[\text{Zeta}\left[\frac{1}{2} + i\, t\right]\right]\right] \text{Sech}[\pi\, t]$$

```
Plot[g, {t, 0, 1.9}]
```

```
NIntegrate[g, {t, 0, 200}, MaxRecursion -> 350, WorkingPrecision -> 200] // Timing
```

NIntegrate::slwcon :
  Numerical integration converging too slowly; suspect one of the following: singularity, value of the
    integration is 0, highly oscillatory integrand, or WorkingPrecision too small. ≫

NIntegrate::eincr :
  The global error of the strategy GlobalAdaptive has increased more than 400 times. The global error is expected
    to decrease monotonically after a number of integrand evaluations. Suspect one
    of the following: the working precision is insufficient for the specified precision
    goal; the integrand is highly oscillatory or it is not a (piecewise) smooth function;
    or the true value of the integral is 0. Increasing the value of the GlobalAdaptive
    option MaxErrorIncreases might lead to a convergent numerical integration.
    NIntegrate obtained ≪258≫ and ≪261≫ for the integral and error estimates. ≫

{48347.532`,
 0.08346229122167157875281209580070247873696869522728842200951359415806355561602135916 \
 8998294841727586222632551828572449564159485875107946103249188344466128712708605660247 \
 9816806612721149070475152531389797126399041487200. }

{48347.5,
 0.08346229122167157875281209580070247873696869522728842200951359415806355561602135916 \
 8998294841727586222632551828572449564159485875107946103249188344466128712708605660247 \
 9816806612721149070475152531389797 }

```
%[[2]] - N[Log[Pi/2] + Sum[(-1)^{n+1} Log[Zeta[n]], {n, 2, 1000}], 200]
```

{1.1626085474678047184389300550686152480014451792021891498879166540003262527669894586 \
  67874108336290709909019813393866151810^-81}

Fig. 1. Computer transcript illustrating the numerical testing of eq. (4).